\documentclass[12pt]{amsart}

\usepackage{amsmath}
\usepackage{amsfonts}
\usepackage[utf8]{inputenc}
\usepackage[T1]{fontenc}
\usepackage{amsthm}
\usepackage{tikz-cd}
\usepackage{amssymb}

\newtheorem{lemma}{Lemma}[section]

\newtheorem{thm}{Theorem}[section]
\newtheorem{prop}{Proposition}[section]

\newtheorem{cor}{Corollary}[section]

\theoremstyle{remark}
\newtheorem{remark}{Remark}[section]

\theoremstyle{claim}

\newcommand{\OO}{\mathcal{O}}
\newcommand{\HH}{\mathcal{H}}
\newcommand{\F}{\mathcal{F}}
\newcommand{\G}{\mathcal{G}}
\newcommand{\U}{\mathcal{U}}
\newcommand{\V}{\mathcal{V}}
\newcommand{\KK}{\mathcal{K}}
\newcommand{\im}{\mathrm{im}\,}
\newcommand{\LIM}{\lim_{\underset{\alpha}{\longrightarrow}}}
\newcommand{\hdim}{\mathrm{hdim}\,}

\newcommand{\pd}{\mathrm{pd}\,}

 \title[]{Blown-up \v{C}ech cohomology and Cartan's Theorem B on real algebraic varieties}
 \author{Tomasz Kowalczyk}
 \date{}
\begin{document}
\maketitle
\begin{abstract}
We introduce a concept of blown-up \v{C}ech cohomology for coherent sheaves of homological dimension $\leq 1$ and some quasi-coherent sheaves on a non-singular real affine variety. Its construction involves a directed set of multi-blowups. We establish, in particular, long exact cohomology sequence and Cartan's Theorem B. Finally, some applications are provided, including universal solution to the first Cousin problem (after blowing up).

\end{abstract}
\section{Introduction}
It is known that Cartan's Theorems A and B do not hold in the real algebraic case. Nevertheless, in \cite{TK} we established a real version of Theorem A to the effect that for any coherent sheaf $\F$ on a non-singular real affine variety $X$ there exists a multi-blowup $\sigma : X_\sigma \rightarrow X$ such that the pullback sheaf $\sigma^*\F$ is generated by global sections on $X_\sigma$.  

In this paper we intend to develop a theory of the so-called blown-up \v{C}ech cohomology groups of quasi-coherent sheaves of global presentation and with stalks of projective dimension $\leq 1$. In particular, we provide long exact cohomology sequence and a real version of Theorem B.

\begin{remark}
There are several cases where Theorem B holds: in complex analytic geometry, in algebraic geometry over algebraically closed field (Serre \cite{Serre}), in scheme theory (Grothendieck \cite{Gr}, Hartshorne \cite{Hartshorne}) as well as recent versions in real regulous geometry (Fichou--Huismann--Mangolte--Monnier \cite{4 x fr}) and in regulous geometry over Henselian valued fields (Nowak, \cite{KN} \cite{KN2}). Note also that the theory of regulous functions is closely related to that of continuous hereditarily rational functions, developed by Koll\'ar--Nowak \cite{K+N}. They used blowups to improve the functions under study, which is also a basic tool in this paper.
\end{remark}

In Section 2, we make the set of all multi-blowups of $X$ a directed set. It is crucial for defining the aforementioned cohomology groups, which are the direct limits with respect to that directed set. In Section 3, we prove that the section functor is right exact after blowing up. In Section 4, we provide basic properties of homological dimension of a sheaf. Section 5 develops the concept of blown-up \v{C}ech cohomology. Section 6 provides an application to the first Cousin problem; more precisely, we show that it is universally solvable after blowing up. Section 7 contains an example of quasi-coherent locally free sheaf $\F$ of infinite rank on $\mathbb{R}^2$, such that after any blowup $\alpha : X_\alpha \rightarrow \mathbb{R}^2$, the pull-back $\alpha^*\F$ cannot be generated by global sections.

Our approach combines  the technique of coherent algebraic sheaves and their \v{C}ech cohomology, developed by Serre \cite{Serre} on algebraic varieties over algebraically closed fields, and transformation to a normal crossing by blowing up.

Throughout the paper, $X$ will be a quasi-projective (hence affine) real algebraic variety with the structure sheaf $\OO_X$ of regular functions. We remind the reader that a sheaf $\F$ of $\OO_X$-modules is quasi-coherent if there exists an open covering $\{U_i \}_{i=1}^n$ of $X$ such that for each $i$ there is an exact sequence of sheaves:
$$ \OO_X^{\oplus J_i}|_{U_i} \xrightarrow{\phi_i} \OO_X^{\oplus I_i}|_{U_i} \xrightarrow{\psi_i} \F|_{U_i} \rightarrow 0.$$

It is clear that one can take a common $J$ with the biggest cardinality among $J_i$ for $i=1,2, \dots n$. A sheaf $\F$ is coherent if $J$ and each $I_i$ can be taken finite.
 A sequence of sheaves is exact if it is exact on stalks. If $\sigma : Y \rightarrow X$ is a morphism of varieties, by $\sigma^*s$ or $(s)^\sigma$ we denote the pull-back of a section $s$ or a function $f$, respectively. For $U \subset X$, $U^\sigma:= \sigma^{-1}(U)$ denotes the preimage under $\sigma$ of $U$; similarly, for an open covering $\U = \{U_i \}^n_{i=1}$, we put $\U^\sigma = \{U_i^\sigma \}_{i=1}^n$.

From now on we shall assume that a given real affine variety $X$ is irreducible and non-singular. In case of a non-singular reducible variety, every reasoning can be carried out on each component separately. By a multi-blowup $\sigma : X_\sigma \rightarrow X$ we mean a finite composition of blowups along smooth centres. We say that a regular function $g:X \rightarrow \mathbb{R}$ on a non-singular real algebraic variety of dimension $d$ is a simple normal crossing if, in a neighbourhood of each point $a \in X$, one has 
$$  g(x)=u(x)x^\alpha=u(x)x_1^{\alpha_1}x_2^{\alpha_2}\dots x_d^{\alpha_d},   $$
where $u(x)$ is a unit at $a$, $\alpha \in \mathbb{N}^d$ and $x=(x_1,x_2,\dots,x_d)$ are local coordinates near $a$, i.e.\ $x_1,x_2,\dots, x_d \in \OO_{a,X}$ is a regular system of parameters of the local ring  $ \OO_{a,X}$.

One of the  basic tools applied in this paper is transformation to a simple normal crossing by blowing up recalled below (see e.g.\ \cite[Theorem 3.26]{Kollar}).

\begin{thm}\label{Thm 1}
Let $f : X \rightarrow \mathbb{R}$ be a regular function on a non-singular real algebraic variety $X$. Then there exists a multi-blowup $\sigma : X_\sigma \rightarrow X$ such that $f^\sigma:=f \circ \sigma$ is a simple normal crossing.
\hfill$\square$
\end{thm}

A useful strengthening of this theorem, stated below, relies on the following elementary result (see e.g.\ \cite[Lemma 4.7]{bierstone}.

\begin{lemma} \cite[Lemma 4.7]{bierstone}
Let $x=(x_1,x_2,\dots, x_p)$ be a regular system of parameters of $\OO_{x,X}$. Let $\alpha, \beta,\gamma \in \mathbb{N}^p$ and let $a(x),b(x),c(x)$ be invertible elements in $\OO_{x,X}$. If
$$ a(x)x^\alpha - b(x)x^\beta = c(x)x^\gamma,$$
then either $\alpha \leq \beta$ or $\beta \leq \alpha$. Here inequality $\alpha \leq \beta$ means that $\alpha_j \leq \beta_j$ for all $j=1,2,\dots, p$.
\hfill$\square$
\end{lemma}

We immediately obtain

\begin{cor}\label{cor 1}
Let $f_1,f_2,\dots, f_k$ be regular functions on $X$. Then there exists a multi-blowup $\sigma : X_\sigma \rightarrow X$ such that $f^\sigma_1,f^\sigma_2,\dots,f^\sigma_k$ simultaneously are simple normal crossings which locally are linearly ordered by divisibility relation near each point $b \in X_\sigma$.
\begin{proof}
Apply Theorem \ref{Thm 1} to the function 
$$f:=f_1f_2\dots f_k \prod_{i<j}(f_i-f_j).$$
Then all of the functions $f^\sigma_j$ and $f^\sigma_i - f^\sigma_j, \ i,j=1,2, \dots,k, i<j,$ simultaneously are simple normal crossing. Now the conclusion follows directly from the above lemma.

\end{proof}
\end{cor}
Throughout the paper, greek letters $\alpha, \beta, \gamma, \sigma$ will denote multi-blowups of $X$, $X_\alpha, X_\beta, X_\gamma, X_\sigma$ their domains and calligraphic letters $\F, \G, \HH$ sheaves on $X$.

Note finally that blowups of real affine varieties remain affine, as so are real projective spaces (\cite[Theorem 3.4.4]{Bochnak}).

\section{Directed system of multi-blowups}

Given two multi-blowups $\alpha : X_\alpha \rightarrow X, \ \beta : X_\beta \rightarrow X$ we say that $ X_\alpha \succeq X_\beta $ if there is a (unique) regular map $f_{\alpha \beta}$ making the following diagram commute
$$
\begin{tikzcd}
 X_\alpha \arrow[rd, "\alpha"] \arrow[r, "f_{\alpha \beta}"] & X_\beta  \arrow[d, "\beta"] \\
& X
\end{tikzcd}
$$
Obviously, $\succeq$ is a reflexive and transitive relation on the set of all multi-blowups of $X$.
\begin{prop}\label{Prop 1}
With the relation given above, the set of multi-blowups of $X$ is a directed set.
\begin{proof}
We need to show that for any two multi-blowups $\sigma_1 :X_1 \rightarrow X$ and $ \sigma_2 :X_2 \rightarrow X$ there is a multi-blowup $\sigma_3 :X_3 \rightarrow X$ such that $X_3 \succeq X_1$ and $X_3 \succeq X_2$.
Let $\phi : X_1 \dashrightarrow X_2$ be a rational map which makes the diagram
$$
\begin{tikzcd}
& X_1 \arrow[r, dashed, "\phi"] \arrow[rd, "\sigma_1"] & X_2  \arrow[d, "\sigma_2"]\\
& & X \\
\end{tikzcd}
$$
commutative. Let $\mathrm{dom}(\phi)$ be the biggest Zariski open subset of $X_1$ on which $\phi$ is a regular mapping. Clearly, $X_2$ is an affine variety embedded into $\mathbb{R}^N\subset \mathbb{P}^N(\mathbb{R})$ for some $N$; embed $\mathbb{R}^N $ into $\mathbb{P}^N(\mathbb{R}) $ by the map
$$(x_1,x_2,\dots,x_N) \mapsto (x_1:x_2:\dots:x_N:1).$$
Then $\phi$ can be treated as a map into $\mathbb{R}^N$, $\phi : X_1 \dashrightarrow \mathbb{R}^N$,
 with a presentation
  $$\phi(x)=\left( \frac{\phi_1(x)}{q(x)}, \frac{\phi_2(x)}{q(x)}, \dots, \frac{\phi_N(x)}{q(x)}\right),$$
   where $\phi_1, \dots, \phi_N,q$ are regular on $X_1$ and $\{q=0 \} \cap \mathrm{dom}(\phi)= \emptyset$. Consider a multi-blowup $\tau : X_3 \rightarrow X_1$ from Corollary \ref{cor 1} applied to the regular functions $\phi_1, \dots, \phi_N,q$ on $X_1$. Then we get the commutative diagram
$$
\begin{tikzcd}
& & X_3 \arrow[dl, "\tau" ] \arrow[d, dashed, "\phi \circ \tau"] \arrow[rdr, "f"] \\
 & X_1 \arrow[r, dashed, "\phi"] \arrow[rd, "\sigma_1"]  & X_2 \arrow[r, hook] \arrow[d, "\sigma_2"] & \mathbb{R}^N \arrow[r, hook,] & \mathbb{P}^N & \\
   &  & X
\end{tikzcd}
$$
with the function $f : X_3 \rightarrow \mathbb{P}^N(\mathbb{R})$ given by the formula
$$f = (\phi_1 \circ \tau : \dots : \phi_N \circ \tau : q \circ \tau).$$
It follows immediately from the conclusion of Corollary \ref{cor 1} that the map $f$ is regular. Note that the maps $\sigma_1, \sigma_2, \tau$ are proper.

We must show that $f(X_3)\subset X_2$. If $y \in X_3$, $\tau(y) \in \mathrm{dom}(\phi)$, then $f(y) \in X_2$ by the commutativity of the upper left triangle diagram. Assume that $\tau(y) \notin \mathrm{dom}(\phi)$. Let $\{y_n\} \subset X_3$ be a sequence converging to $y$ in the Euclidean topology such that $\tau(y_n)\in \mathrm{dom}(\phi)$. Such a sequence can always be found because the blowups are biregular on a Zariski open subset and a Zariski open subsets of a non-singular irreducible variety are dense in the Euclidean topology. Let  
$$K:=\{(\sigma_1 \circ \tau)(y_n) : n \in \mathbb{N}\} \cup  \{(\sigma_1 \circ \tau)(y) \} \subset X$$
and
$$ V := \sigma_2^{-1}(K) \subset X_2.$$
Obviously $K$ is a compact set, so is $V$ as $\sigma_2$ is proper. We get 
$$ \sigma_2(f(y_n))=\sigma_2(\phi(\tau(y_n)))=\sigma_1(\tau(y_n)) \in K.$$
Hence by the definition of $V$, $f(y_n) \in V$. Thus the sequence $\{f(y_n) \}$ is convergent to $f(y) \in V \subset X_2$, as desired.

\end{proof}
\end{prop}
Proposition \ref{Prop 1} will be crucial for the construction of blown-up \v{C}ech cohomology given in section 5.
\begin{remark}
The above construction for proper schemes was described in \cite{shafarevich}. It was broadly used in the theory of real holomorphy rings (see e.g. \cite{buchner kucharz}).
\end{remark}

\section{Section functor is right exact after blowing up}

Let $\tau : Y \rightarrow X$ be a morphism of varieties and let $\F$ be a sheaf of $\OO_X$-modules. For any Zariski open set $U\subset X$ we have the canonical and functorial homomorphism
$$ \tau^*: \F(U) \rightarrow \tau^*\F(U^\tau)$$ 
such that
$$(\tau^*s)(y)=s({\tau(y)}) \otimes 1 \in (\tau^{-1}\F)_{\tau(y)} \otimes_{\OO_{\tau(y),X}} \OO_{y, Y}=(\tau^*\F)_y, $$
for $s \in \F(U)$ and each $y \in U^\tau$. The main aim of this section is to prove

\begin{thm}\label{Thm 3}
Let 
$$ \G \overset{\theta}{\rightarrow} \HH \rightarrow 0$$ 
be an exact sequence of quasi-coherent sheaves on $X$. Then for any Zariski open subset $U\subset X$ and any section $u \in \HH(U)$ there exists a multi-blowup $\alpha : X_\alpha \rightarrow X$ such that $\alpha^*u \in \im  \theta^\alpha$.
\end{thm}
 
A crucial role in the proof of Theorem \ref{Thm 3} is played by Lemmas \ref{q-coh ext zero section} and \ref{q-coh ext section to global}. Before proving them, we recall three Corollaries 2.1, 2.2, 2.3 from our paper \cite{TK}. The first two will be used in the proofs of Lemmas \ref{q-coh ext zero section} and \ref{q-coh ext section to global}, and the last (being a corollary to Cartan's Theorem A) in Section 5.

\begin{lemma}\cite[Corollary 2.1]{TK}\label{lokalne przedluzanie funkcji regularnych}
Let $Q$ be a regular function on $X$. Then for any finite number of regular functions $P_i \in \OO_X(X)$ such that $V(P_i)\subset V(Q)$ for $i=1,2,\dots,k$,  there exist a multi-blowup $\sigma : X_\sigma \rightarrow X$ and a positive integer $N$ such that 
$$ (Q^N)^\sigma \in P_i^\sigma \cdot \OO_{X_\sigma}(X_\sigma)$$
for each $i=1,2,\dots, k$.
\hfill$\square$
\end{lemma}

\begin{lemma} \cite[Corollary 2.2]{TK}\label{przedluzanie funkcji regularnych}
Let $Q$ be a regular function on $X$ and $U=X\setminus \{Q=0 \}$. Then for any $f \in \OO_X(U)$ there exist a multi-blowup $\sigma : X_\sigma \rightarrow X$ and a positive integer $N$ such that $(Q^Nf)^\sigma$  can be extended to a global regular function on $X_\sigma$.
\hfill$\square$
\end{lemma}

\begin{lemma} \cite[Corollary 2.3]{TK}\label{Lemma 8}
Let $\F$ be a coherent sheaf on $X$. Then there exists a multi-blowup $\sigma : X_\sigma \rightarrow X$ such that the pull-back $\sigma^*\F$ admits a global presentation:
$$ \OO^p_{X_\sigma} \rightarrow \OO^q_{X_\sigma} \rightarrow \sigma^*\F \rightarrow 0.  $$
\hfill$\square$
\end{lemma}

\begin{lemma}\label{q-coh ext zero section}
Let $\F$ be a quasi-coherent sheaf on $X$. For any $Q \in \OO_X(X)$ and a section $s \in \F(X)$ such that $s|_{U} =0$ with $U=X\setminus\{ Q =0 \}$, there exist a multi-blowup $\sigma : X_\sigma \rightarrow X$ and a positive integer $N$ such that $(Q^N)^\sigma \sigma^*s = 0$ in $\sigma ^* \F (X_\sigma)$.

\begin{proof}
By quasi-compactness, there is a finite Zariski open covering $\{U_i \}_{i=1}^n$ of $X$ such that for each $i$ we have:
\begin{itemize}
\item[1)]
a presentation
$$ \OO_X^{\oplus J}|_{U_i} \xrightarrow{\phi_i} \OO_X^{\oplus I_i}|_{U_i} \xrightarrow{\psi_i} \F|_{U_i} \rightarrow 0,$$
\item[2)]
$s|_{U_i}=\psi_i(t_i)$ for some $t_i \in \OO_X^{\oplus I_i}(U_i).$
\end{itemize}
Put 
$$ \mathrm{Rel}\,(t_i,\phi_i(e_j): j \in J; \OO_X(U_i)):=$$ 
$$=\left\{(q,(q_j)_{j\in J}) \in \OO_X(U_i) \oplus \OO_X(U_i)^{\oplus J}\ : qt_i +  \sum_{j \in J} q_j \phi_i(e_j)=0  \right\}= $$
$$=\lim_{\substack{\longrightarrow \\ J_f \subset J \\ \# J_f <\infty}} \mathrm{Rel}\,(t_i,\phi_i(e_j): j \in J_f; \OO_X(U_i)).$$
In other words, we express $\mathrm{Rel}\,(t_i,\phi_i(e_j): j \in J; \OO_X(U_i))$ as a direct limit of modules of relations of $t_i$ and finitely many elements indexed by $J$. Let
$$\mathcal{I}_i:= \left\{ q \in \OO_X(U_i) : \exists (q_j)_{j \in J_f}, \# J_f < \infty:qt_i + \sum_{j \in J_f} q_j \phi_i(e_j)=0 \right\}= $$
$$={\pi_1 \left(\mathrm{Rel}\,(t_i,\phi_i(e_j): j \in J; \OO_X(U_i))\right) \subset \OO_X(U_i)}$$
for each $i=1,2,\dots,n$; here $\pi_1$ is the natural projection onto the first factor, $e_j$ is an element of $\OO_X(U_i)^{\oplus J}$ which has $1$ on $j$-th entry and zero elsewhere.
Then, for every $x \in U_i$ we have
$$ \mathrm{Rel} \, (t_i(x), \phi_i(e_j(x), j \in J; \OO_{x,X})=$$
$$=\lim_{\substack{\longrightarrow \\ J_f \subset J \\ \# J_f <\infty}} \mathrm{Rel} \, (t_i(x), \phi_i(e_j)(x), j \in J; \OO_{x,X})= $$
$$= \lim_{\substack{\longrightarrow \\ J_f \subset J \\ \# J_f <\infty}} \Big( \mathrm{Rel}\,(t_i,\phi_i(e_j): j \in J_f; \OO_X(U_i)) \cdot \OO_{x,X}  \Big)=$$
$$=\Big(\lim_{\substack{\longrightarrow \\ J_f \subset J \\ \# J_f <\infty}}  \mathrm{Rel}\,(t_i,\phi_i(e_j): j \in J_f; \OO_X(U_i))\Big) \cdot \OO_{x,X} =   $$
$$= \Big( \mathrm{Rel}\,(t_i,\phi_i(e_j): j \in J; \OO_X(U_i)) \Big) \cdot \OO_{x,X}  $$
because modules of relations commute with flat base change (see e.g.\ \cite{Bourbaki}, Chap.\ I, \S2, Remark 2 after Prop.\ 1).
Therefore,
$$ \mathcal{I}_i \cdot \OO_{x,X}= \{ q_x \in \OO_{x,X} : \exists (q_{jx})_{j \in J_f}, \# J_f <\infty , qt_i(x) + \sum_{j \in J_f}  q_j\phi_i(e_j)(x)=0 \}$$
and thus, $1 \in \mathcal{I}_i \cdot \OO_{x,X}$ for every $x \in U_i \cap U$. Hence we get
$$U_i\cap V(\mathcal{I}_i) \subset U_i \cap V(Q).$$

Clearly, there exists $p_i \in I_i$ such that $V(p_i)=V(\mathcal{I}_i).$ Each $p_i$ may be written in the form
$$  p_i = \frac{P_i}{R_i} \,\,\, \mathrm{with} \,\,\, P_i,R_i \in \OO_X(X) \,\,\, \mathrm{and} \,\,\, V(R_i)\cap U_i = \emptyset$$
for $i=1,2,\dots, n$. Hence
$$V(p_i)\cap U_i=V(P_i)\cap U_i$$
and by Lemma \ref{lokalne przedluzanie funkcji regularnych} there exist a multi-blowup 
$\sigma : X_\sigma \rightarrow X$ and a positive integer $N$ such that
$$ (Q^N)^\sigma|_{U_i^\sigma} \in P_i^\sigma|_{U_i^\sigma} \cdot \OO_{X_\sigma}(U_i^\sigma), \ \ i=1,\dots,n.$$
As $p_i$ and $P_i$ differ only by a unit on $U_i$, we get 
$$ (Q^N)^\sigma|_{U_i^\sigma} \in p_i^\sigma|_{U_i^\sigma} \cdot \OO_{X_\sigma}(U_i^\sigma)\ \ i=1,\dots, n.$$
Therefore
\begin{equation*}
\begin{split}
 (Q^N)^\sigma|_{U_i^\sigma} \cdot \sigma^*(t_i) \in p_i^\sigma|_{U_i^\sigma} \cdot
 \sigma^*(t_i) \cdot \OO_{X_\sigma}(U_i^\sigma)= \\
 =\sigma^*(p_i|_{U_i}t_i)\cdot \OO_{X_\sigma}(U_i^\sigma)
  \subset \sigma^*(\phi_i(\OO_X^{\oplus J}(U_i)))=
  \phi_i^\sigma(\OO_{X_\sigma}^{\oplus J}(U_i^\sigma)), 
\end{split}
\end{equation*}
whence, 
$$(Q^N)^\sigma|_{U_i^\sigma}\cdot \sigma^*s|_{U_i^\sigma}=\psi_i^\sigma((Q^N)^\sigma\sigma^*t_i) = 0|_{U_i^\sigma}.$$
This finishes the proof.

\end{proof}
\end{lemma}

\begin{lemma}\label{q-coh ext section to global}
Let $\F$ be a quasi-coherent sheaf on $X$ with local presentations
$$ \OO_X^{\oplus J}|_{U_i} \xrightarrow{\phi_i} \OO_X^{\oplus I_i}|_{U_i} \xrightarrow{\psi_i} \F|_{U_i} \rightarrow 0 \ \ i=1,2,\dots,n$$
 on a finite Zariski open covering $\{U_i \}^n_{i=1}$ of $X$.
 Consider a finite number of sections $s_j \in \F(V_j)$ on Zariski open sets
  $$V_j=X\setminus\{Q_j=0\}, \,\,\,\, j=1,2,\dots,m$$
  where $Q_j$ are regular functions on $X$. Assume that every $V_j$ is contained in $U_{i(j)}$ for some $i(j)=1,2,\dots, n$ and  that for each $j$ there is a section $t_j \in \mathcal{O}_X^{\oplus I_{i(j)}}(V_j)$ such that $\psi_{i(j)}(t_j)=s_j$. Then there exists a positive integer $N$ and a multi-blowup $\sigma : X_\sigma \rightarrow X$ such that  every section $(Q_j^N)^\sigma \sigma^*s_j, \ \  j=1,2,\dots, m$, extends to a global section on $X_\sigma.$

\begin{proof}
Since taking pull-back under a multi-blowup does not affect the assumptions, it suffices to consider only one $j=1,2,\dots,m$. So fix an index $j$, and let $t_{ji}=t_j|_{U_i\cap V_j}=(t_{jil})_{l \in I_i}$. Since $t_{ji}$ is an element of direct sum, all but finitely many $t_{jil}$ are zero, so every $t_{ji}$ can be identified with finitely many regular functions on $U_i \cap V_j$, indexed by some subset $J_i \subset I_i$ such that $\# J_i < \infty$ i.e.\ $t_{ji}=(t_{jil})_{l \in J_i}$.
We have 
$$t_{jil}=\frac{t_{jil1}}{t_{jil2}}, \ t_{jil1}, t_{jil2} \in \OO_X(X)$$
and
$$V(t_{jil2})\cap U_i \subset V(Q_j)\cap U_i.$$
Using Lemma \ref{przedluzanie funkcji regularnych} we can find a positive integer $N_1$ and a multi-blowup $\alpha: X_\alpha \rightarrow X$ such that
$$((Q_j)^{N_1}t_{jil})^{\alpha} \in \OO_{X_1}({U_i^{\alpha}}) \ \ \text{for all} \ \ j,i,l.$$
 Now define $\widetilde{s_{ji}}:=\psi^{\alpha}_{i}(((Q_j)^{N_1}t_{ji})^{\alpha})$. Then for any two distinct indices $i_0, \, i_1 \in \{1,2,\dots,n\} $ we have
$$ (\widetilde{s_{ji_0}}-\widetilde{s_{ji_1}})|_{U^{\alpha}_{i_0} \cap U^{\alpha}_{i_1}}  \in (\alpha^*\F)(U^{\alpha}_{i_0} \cap U^{\alpha}_{i_1})$$
and
$$ (\widetilde{s_{ji_0}}-\widetilde{s_{ji_1}})|_{U^{\alpha}_{i_0} \cap U^{\alpha}_{i_1} \cap V^{\alpha}_j} =0.   $$
By Lemma \ref{q-coh ext zero section}, we can find a multi-blowup $\beta : X_\beta \rightarrow X_\alpha$ and a positive integer $N_2$ such that
$$( (Q_j^{N_1 + N_2})^{\alpha \circ \beta}\beta^*\widetilde{s_{ji_0}}  - (Q_j^{N_1 + N_2})^{\alpha \circ \beta}\beta^*\widetilde{s_{ji_1}})|_{U^{\alpha \circ \beta}_{i_0} \cap U^{\alpha \circ \beta}_{i_1} } = 0. $$

Considering all distinct pairs of indices $i_0,i_1$, we can assume that the differences as above vanish for all those pairs. Therefore the sections 
$$((Q_j^{N_1 + N_2})^{\alpha \circ \beta}\beta^*\widetilde{s_{ji}} )|_{{U^{\alpha \circ \beta}_i}}, \ \ i=1,2,\dots,n, $$
glue together to a global section on $X_\beta$. Thus 
$$\sigma:=\alpha \circ \beta :X_\beta=X_\sigma\rightarrow X$$
 is the multi-blowup we are looking for.

\end{proof}

\end{lemma}

\noindent
{\em Proof of Theorem 3.1}
Let $\{ V_i\}_{i=1}^n$ be a finite Zariski open covering of $U$ satisfying following two conditions 
\begin{itemize}
\item[a)]
$$ \OO_X^{\oplus J}|_{V_i} \xrightarrow{\phi_i} \OO_X^{\oplus I_i}|_{V_i} \xrightarrow{\psi_i} \G|_{V_i} \rightarrow 0 \ \ i=1,2,\dots,n$$

\item[b)]
there exist sections $s_i \in \G(V_i)$ and $t_i \in \OO_X^{\oplus I_i}(V_i)$ such that $\theta(s_i)=u|_{V_i}$ and $\psi_i(t_i)=s_i$.

\end{itemize}

 Each $V_j$ is of the form $V_j=U\setminus \{Q_j=0 \}$ for some $Q_j \in \OO_X(X)$.

By Lemma \ref{q-coh ext section to global}, there exist a multi-blowup $\alpha : X_\alpha \rightarrow X$ and a positive integer $N$ such that every section $(Q_j^N)^\alpha \alpha^*s_j$ extends to $U^\alpha$:
$$(Q_j^N)^\alpha \alpha^*s_j \in \alpha^*\G(U_\alpha).$$
Obviously,
$$ \left(\theta^\alpha((Q_j^N)^\alpha \alpha^*s_j)-(Q_j^N)^\alpha \alpha^*u\right)|_{V_j^\alpha}=0$$
for each $j=1,2,\dots,m$. By Lemma \ref{q-coh ext zero section} there exists a multi-blowup $\beta : X_\beta \rightarrow X_\alpha$ and a positive integer $M$ such that
$$\theta^{\alpha \circ \beta}\left((Q_j^{N+M})^{\alpha \circ \beta} (\alpha \circ \beta)^*s_j\right)=(Q_j^{N+M})^{\alpha \circ \beta} (\alpha \circ \beta)^*u$$
for each $j=1,2, \dots,m$. We may, of course, assume that $N+M$ is even. Then the function 
$$ R:= \frac{1}{\sum_{j=1}^m (Q_j^{N+M})^{\alpha \circ \beta}} $$
is regular on $U^{\alpha \circ \beta}$ , because $V_j$ was a covering of $U$.
We have
$$ R \sum_{j=1}^m \theta^{\alpha \circ \beta}((Q_j^{N+M})^{\alpha \circ \beta} (\alpha \circ \beta)^*s_j) =R\sum_{j=1}^m (Q_j^{N+M})^{\alpha \circ \beta} (\alpha \circ \beta)^*u = (\alpha \circ \beta)^*u. $$
This concludes the proof.

\hspace*{\fill} $\square$
\newline

Now we state two direct consequences of Theorem \ref{Thm 3}.

\begin{cor}

Let 
$$ \G \overset{\theta}{\rightarrow} \HH \rightarrow 0$$ 
be an exact sequence of quasi-coherent sheaves on $X$. Consider $\{U_i \}_{i=1}^m$ a finite collection of Zariski open subsets and sections $u_i \in \HH(U_i)$ for $i=1,2,\dots,m$. Then there exists a multi-blowup $\alpha : X_\alpha \rightarrow X$ such that for each $i, \ \alpha^*u_i \in \im  \theta^\alpha$.
\begin{proof}
This can be obtained by repeated application of Theorem \ref{Thm 3}.
\end{proof}
\end{cor}

Let $\F$ be a quasi-coherent sheaf on $X$. Put
$$\widetilde{\F}(U) = \LIM \alpha^* \F(U^\alpha),$$
for any Zariski open subset $U$ of $X$; direct limit is taken over the directed set of multi-blowups of $X$. Obviously, $\widetilde{\F}(U)$ has a structure of $\OO_X(U)$-module.

It is clear that any element of $\widetilde{\F}(U)$ can be represented as a class $[s]$ for some multi-blowup $\alpha : X_\alpha \rightarrow X$ and some section $s \in \alpha^* \F(U^\alpha)$.

\begin{cor}\label{wniosek z surjekcja na sekcjach}
Let 
$$ \G \overset{\theta}{\rightarrow} \HH \rightarrow 0$$ 
be an exact sequence of quasi-coherent sheaves on $X$. Then, for any Zariski open subset $U\subset X$, the induced sequence of $\OO_X(U)$-modules
$$\widetilde{\G}(U) \rightarrow \widetilde{\HH}(U) \rightarrow 0 $$
is exact.
\begin{proof}
The above exact sequence is well defined since $\alpha^*$ is functorial.
Let $[s] \in \widetilde{\HH}(U)$ with a representative $ s \in \alpha^*\HH(U^\alpha)$ for some multi-blowup $\alpha : X_\alpha \rightarrow X$. To finish the proof, it is enough to apply Theorem \ref{Thm 3} to $s$.

\end{proof}
\end{cor}

\section{Homological dimension of sheaves}
In this section, we assume all rings to be commutative with unity. By a local ring we mean a ring with unique maximal ideal. Let $\HH$ be a quasi-coherent sheaf on $X$ and $x \in X$.

We say that $\HH$ is of \emph{homological dimension $k$ at $x$}, $\hdim_x\HH=k$, if $k$ is the smallest integer such that there exists a Zariski open neighbourhood $U$ and sets of indices $I_0,I_1,\dots,I_k$ for which there is an exact sequence of sheaves:
$$ 0 \rightarrow \OO_X^{\oplus I_k}|_U \rightarrow \OO_X^{\oplus I_{k-1}}|_U \rightarrow \dots \rightarrow \OO_X^{\oplus I_0}|_U \rightarrow \HH|_U \rightarrow 0.$$

We define the \emph{homological dimension of $\HH$} as 
$$\hdim \HH=\sup_{x \in X} \hdim_x \HH.$$

 
Obviously, $\hdim \HH = 0$ iff $\HH$ is a locally free sheaf. Consequently, $\hdim \HH = 1$ means that $\HH$ is locally a quotient of free sheaves.

\begin{remark}\label{remark 0}
In the case of a coherent sheaf $\HH$, all sets $I_j$ can be finite.
\end{remark}

\begin{remark}\label{remark 1}
 If $\HH$ is of homological dimension $k$ at $x$ and $U$ is a sufficiently small open neighbourhood of $x$, then $\hdim_y\HH\leq k$ for $y \in U$. Then for any point $y \in U$ the $\OO_{y,X}$-module $\HH_y$ is of projective dimension $\pd \HH_y  \leq k$, because the notions of free and projective modules coincide over a local ring (cf.\ \cite{kaplansky 2}).
\end{remark}

\begin{prop}\label{tw z wymiar hdim=pdim}
Let $\HH$ be a coherent sheaf on $X$ and $x \in X$, then 
$$\pd\HH_x=\hdim_x \HH,$$
here $ \pd\HH_x$ is a projective dimension of $\OO_{x,X}$-module.

\begin{proof}
The inequality $\leq$ follows from Remark \ref{remark 1}.
Conversely, let $k=\pd\HH_x$. Then there is an exact sequence of $\OO_{x,X}$-modules
$$ 0 \longrightarrow \OO_{x,X}^{p_k} \overset{\phi_{xk}}{\longrightarrow} \dots \longrightarrow \OO_{x,X}^{p_1} \overset{\phi_{x1}}{\longrightarrow} \OO_{x,X}^{p_0} \overset{\phi_{x0}}{\longrightarrow} \HH_x \longrightarrow 0.$$
By \cite[I, \S2, Proposition 5]{Serre} the homomorphisms $\phi_{xj}$ can be lifted to homomorphisms of sheaves on some common Zariski open neighbourhood $V$ of $x$:
$$\OO_X^{p_k}|_V \overset{\phi_k}{\longrightarrow}  \dots \OO_X^{p_1}|_V \overset{\phi_1}{\longrightarrow} \OO_X^{p_0}|_V \overset{\phi_0}{\longrightarrow} \HH|_V.$$
By \cite[I, \S2, Theorem 2]{Serre}, $\im \phi_i$ and $\ker \phi_i$ are coherent for all $i=1,2, \dots, k$ and so are the sheaves $ \ker \phi_i / \im \phi_{i+1}$ for all $i=0,1, \dots, k-1$. By assumption, 
$(\ker \phi_i / \im \phi_{i+1})_x=0$. Therefore the above equality propagates to a common Zariski open neighbourhood $U \subset V$ of $x$. Consequently, the sequence is exact
$$0 \rightarrow \OO_X^{p_k}|_U \overset{\phi_k}{\rightarrow}  \dots \OO_X^{p_1}|_U \overset{\phi_1}{\rightarrow} \OO_X^{p_0}|_U \overset{\phi_0}{\rightarrow} \HH|_U \rightarrow 0.$$
Hence  $ \pd\HH_x \geq \hdim _x \HH $, as desired.

\end{proof}
\end{prop}

\begin{remark}\label{rem 4.3 hdim >pd}
If $\HH$ is a quasi-coherent sheaf we can get only $\pd\HH_x \leq \hdim_x \HH$, as not every morphism on stalks can be lifted to a neighbourhood.
\end{remark}

We still need the following result of homological algebra (cf.\ \cite[Part III, Theorem 2]{kaplansky}).

\begin{thm}\label{kaplansky}
Let
$$ 0 \rightarrow F \rightarrow G \rightarrow H \rightarrow 0$$
be a short exact sequence of $R$-modules over a ring $R$.
\begin{itemize}
\item[a)]
If $\pd G > \pd F$, then $\pd H = \pd G$.
\item[b)]
If $\pd G < \pd F$, then $\pd H = \pd F +1$.
\item[c)]
If $\pd G = \pd F$, then $\pd H \leq \pd G +1$.

\end{itemize}
\end{thm}
\hfill$\square$

We can reformulate the above Theorem into a more condensed version due to P.M. Cohn, namely 
$$ \pd G \leq \max ( \pd F, \pd H)$$
with equality unless $\pd H = \pd F +1$.

As a corollary, we obtain the following proposition, which will be useful further in the paper.

\begin{prop}\label{Prop z 4 podpunktami dla modulow}

Let
$$ 0 \rightarrow F \rightarrow G \rightarrow H \rightarrow 0$$
be a short exact sequence $R$-modules over a ring $R$.
\begin{itemize}
\item[a)]
If $\pd  G \leq 1$ and $\pd H \leq 1$ then $\pd F \leq 1$.
\item[b)]
If $\pd F \leq 1$ and $\pd H \leq 1$ then $\pd G \leq 1$.
\item[c)]
If $\pd F=0$ and $\pd G \leq 1$ then $\pd H \leq 1$.
\item[d)]
If $\pd G=0$ and $\pd H \leq 1$ then $\pd F=0$.
\end{itemize}
\begin{proof}
This follows directly from Theorem \ref{kaplansky}.
\end{proof}
\end{prop}

Given a short exact sequence of coherent sheaves, the corollary below indicates the cases where if two of them are of homological dimension $\leq 1$, so is the third one. Actually, this property may not be preserved only in the case where both $\F$ and $\G$ are of homological dimension $1$.

\begin{cor}\label{Lemma 4.2}
Let 
$$ 0 \rightarrow \F \rightarrow \G \rightarrow \HH \rightarrow 0$$
be a short exact sequence of coherent sheaves on $X$.

\begin{itemize}
\item[a)]
If $\hdim \G \leq 1$ and $\hdim \HH \leq 1$, then $\hdim \F \leq 1$.
\item[b)]
If $\hdim \F \leq 1$ and $\hdim \HH \leq 1$, then $\hdim \G \leq 1$.
\item[c)]
If $\hdim \F=0$ and $\hdim \G \leq 1$, then $\hdim \HH \leq 1$.
\item[d)]
If $\hdim \G =0$ and $\hdim \HH \leq 1$, then $\hdim \F =0$.
\end{itemize}

\begin{proof}
Apply Propositions \ref{tw z wymiar hdim=pdim} and  \ref{Prop z 4 podpunktami dla modulow}.
\end{proof}
\end{cor}

\begin{prop}\label{Prop z iniektywnoscia na modulach}
Let $R$ be an integral domain, and $K$ its field of fractions. Consider an exact sequence of $R$-modules 
$$ 0 \rightarrow F \rightarrow G $$
such that $F$ is free and the projective dimension of $G$ is $\leq 1$.
Then for any ring $S$ such that $R\subset S \subset K$ the induced sequence
$$ 0 \rightarrow F \otimes_R S \rightarrow G \otimes_R S $$
is exact.
\begin{proof}
Clearly, $G=F_2 / E$ and $F=F_1/E$ where $E\subset F_1 \subset F_2$ are $R$-modules and $F_2$, $E$ are free. We have the following short exact sequence of $R$-modules
$$ 0 \rightarrow E \rightarrow F_1 \rightarrow F_1/E \rightarrow 0 $$
where $E$ and $F_1/E$ are free. Since free module is projective, the above exact sequence splits and $F_1$ is also a free module. Therefore the canonical homomorphisms
$$ E \otimes S \rightarrow F_1 \otimes S \rightarrow F_2 \otimes S $$
are injective, because $E \hookrightarrow E \otimes S \hookrightarrow E \otimes K$ and the problem of injectivity of free modules reduces to the one for $K$-vector spaces. Consequently, we can regard $E\otimes S $ and $F_1 \otimes S$ as submodules of $F_2 \otimes S$. Hence the canonical homomorphism
$$(F_1 / E) \otimes S \cong (F_1\otimes S)/(E\otimes S) \rightarrow (F_2/E)\otimes S \cong (F_2\otimes S)/(E\otimes S)$$
is injective as asserted.
\end{proof}
\end{prop}

In the proof of Proposition \ref{Cartan B dla pokrycia}, we shall use the following

\begin{cor}\label{cor - pd F =0 pdG <2}
Let 
$$ 0 \rightarrow \F \rightarrow \G$$
be an exact sequence of quasi-coherent sheaves on $X$ such that $\pd\F_x=0$ and $\pd\G_x\leq 1$ for all $x \in X$.

Then for any multi-blowup $\sigma: X_\sigma \rightarrow X$ the induced sequence
$$0 \rightarrow \sigma^*\F \rightarrow \sigma^* \G$$
is exact.

\begin{proof}
Let $\sigma : X_\sigma \rightarrow X$ be a multi-blowup. Take any $x\in X,\ y\in X_\sigma$ such that $\sigma(y)=x$. We have an exact sequence 
$$0 \rightarrow (\sigma^{-1}\F)_y= \F_x \rightarrow (\sigma^{-1}\G)_y=\G_x.$$
Hence and by Proposition \ref{Prop z iniektywnoscia na modulach}, we get
$$0\rightarrow (\sigma^*\F)_y=\F_x \otimes_{\OO_{x,X}} \OO_{y,X_\sigma} \rightarrow
(\sigma^*\G)_y=\G_x \otimes_{\OO_{x,X}} \OO_{y,X_\sigma};$$
the above sequence is exact because $\OO_{y,X_\sigma}$ is a localization of $\OO_{x, X}$.

\end{proof}
\end{cor}

\begin{lemma}\label{hdim=1 jest zachowane przy pullbacku}
Let $\alpha : X_\alpha \rightarrow X$ be a multi-blowup of $X$ and $\HH$ a quasi-coherent sheaf on $X$. If $\HH$ is of homological dimension $\leq 1 $, so is the pull-back $\alpha^*\HH$.

\begin{proof}
Let $\alpha : X_\alpha \rightarrow X$ be a multi-blowup of $X$, take any $x \in X$ and $U \subset X$ as in the definition of homological dimension. We have a short exact sequence of sheaves
$$ 0 \longrightarrow \OO_X^{\oplus I_1}|_U \longrightarrow \OO_X^{\oplus I_0}|_U \longrightarrow \HH|_U \longrightarrow 0.$$
 By the above corollary and right exactness of pull-back we obtain a short exact sequence of sheaves
 $$ 0 \longrightarrow \OO_{X_\alpha}^{\oplus I_1}|_{U^\alpha} \longrightarrow \OO_{X_\alpha}^{\oplus I_0}|_{U^\alpha} \longrightarrow \alpha^*\HH|_{U^\alpha} \longrightarrow 0.$$
Since $x$ was arbitrary we get $\hdim \alpha^*\HH \leq 1$, as asserted.
\end{proof}

\end{lemma}

Corollary \ref{cor - pd F =0 pdG <2} along with Lemma \ref{hdim=1 jest zachowane przy pullbacku} and Remark \ref{rem 4.3 hdim >pd} yield immediately the following
\begin{cor}\label{wniosek z iniektywnoscia pullbacku dla snopow}
Let 
$$ 0 \rightarrow \F \rightarrow \G$$
be an exact sequence of quasi-coherent sheaves on $X$ such that $\hdim \F =0$ and $\hdim \G \leq 1$. Then for any multi-blowup $\sigma: X_\sigma \rightarrow X$ the induced sequence
$$0 \rightarrow \sigma^*\F \rightarrow \sigma^* \G$$
is exact.
\hfill$\square$
\end{cor}

\section{Blown-Up \v{C}ech Cohomology}

In this section, we introduce the concept of blown-up \v{C}ech cohomology for quasi-coherent sheaves of homological dimension $\leq 1$. Our construction combines the classical one due to Serre \cite{Serre} with direct limit with respect to the directed set of multi-blowups described in Section 2.

Let $\F$ be a sheaf on $X$ and $ \U=\{U_i \}_{i=1}^n$ be a finite Zariski open covering of $X$. Put $U_{i_0\dots i_q}=U_{i_0}\cap U_{i_1}\cap \dots \cap U_{i_q}$ and
$$C^q({\mathcal{U}, \F}):=\prod_{1\leq i_0, i_1 , \dots , i_q \leq n} \F(U_{i_0i_1 \dots i_q}).$$
$C^q({\mathcal{U}, \F})$ is called the abelian group of $q$-cochains. We have a chain complex
$$ \dots \longrightarrow C^{q-1}({\mathcal{U}, \F}) \overset{d^{q-1}}{\longrightarrow} C^q({\mathcal{U}, \F}) \overset{d^q}{\longrightarrow} C^{q+1}({\mathcal{U}, \F}) \longrightarrow \dots $$
where  
$$(d^qf)_{i_0 i_1, \dots i_{q+1}} = \sum_{j=0}^{q+1} (-1)^j f_{i_0 i_1 \dots \widehat{i_j} \dots i_{q+1}}|_{U_{i_0 i_1 \dots i_{q+1}}}$$
for any $f=(f_{{i_0\dots i_q}}) \in C^q(\U, \F)$.
If $\sigma : Y \rightarrow X$ is a morphism of real affine varieties, we get the induced chain complex
$$ \dots \longrightarrow C^{q-1}({\U^\sigma, \sigma^*\F}) \overset{d^{q-1}}{\longrightarrow} C^q({\U^\sigma, \sigma^*\F}) \overset{d^q}{\longrightarrow} C^{q+1}({\U^\sigma, \sigma^*\F}) \longrightarrow \dots. $$ 
and a canonical chain complex homomorphism 
$$\sigma^* : C^\bullet(\U,\F) \rightarrow C^\bullet(\U^\sigma, \sigma^*\F);$$
the canonical homomorphism $\F(U) \rightarrow \sigma^*\F(U^\sigma)$ was described in Section 3. Therefore $\sigma$ induces a homomorphism of cohomology complexes
$$\sigma^* : H^\bullet( \U, \F) \rightarrow H^\bullet(\U^\sigma,\sigma^*\F).$$
Let $\widetilde{C}^\bullet(\U,\F)$ be the chain complex defined by the formula 
$$\widetilde{C}^\bullet(\U,\F)=\LIM C^\bullet(U^\alpha, \alpha^*\F)$$
where the limit is taken over a directed system of multi-blowups of $X$.

 \emph{The $q$-th blown-up \v{C}ech cohomology group of $\F$ with respect to $\U$} $\widetilde{H}^q(\U, \F)$ is the $q$-th cohomology group of the chain complex $\widetilde{C}^\bullet(\U,\F)$. 

\begin{remark}
It is well known that direct limit functor commutes with cohomology functor (see e.g.\ \cite[Theorem 4.14]{balc 1}). Hence 
$$\widetilde{H}^q(\U, \F) = \LIM H^q(\U^\alpha, \alpha^*\F).$$

\end{remark}

We now establish long exact cohomology sequence for quasi-coherent sheaves of homological dimension $\leq 1$, which plays a crucial role in the cohomology theory developed in this paper. The general case of arbitrary quasi-coherent sheaves is not at our disposal, as pull-back functor $ \F \mapsto \alpha^*\F$ (along with tensor product functor) is not left exact. But we have of course the following

\begin{prop}

Let $\U$ be a finite Zariski open covering of $X$ and
$$0 \rightarrow \F \rightarrow \mathcal{G} \rightarrow \mathcal{H} \rightarrow 0$$
be a short exact sequence of quasi-coherent sheaves on $X$. Suppose that for any multi-blowup $\sigma : X_\sigma \rightarrow X$ the induced sequence

$$0 \rightarrow \sigma^*\F \rightarrow \sigma^*\mathcal{G}$$ 
is exact. Then there is a short exact sequence of chain complexes

$$ 0 \rightarrow   \widetilde{C}^\bullet({\U, \mathcal{F}}) \rightarrow \widetilde{C}^\bullet({\U, \mathcal{G}}) \rightarrow \widetilde{C}^\bullet({\U, \mathcal{H}}) \rightarrow 0$$
which induces a long exact sequence of blown-up \v{C}ech cohomology with respect to $\U$

$$ \dots \rightarrow \widetilde{H}^p(\U, \F) \rightarrow \widetilde{H}^p(\U, \G) \rightarrow \widetilde{H}^p(\U, \HH) \rightarrow \widetilde{H}^{p+1}(\U, \F) \rightarrow $$
$$ \widetilde{H}^{p+1}(\U, \G) \rightarrow \widetilde{H}^{p+1}(\U, \HH) \rightarrow \widetilde{H}^{p+2}(\U, \F) \rightarrow \dots $$

\begin{proof}
Of course the induced short sequence

$$0 \rightarrow \sigma^*\F \rightarrow \sigma^*\mathcal{G} \rightarrow \sigma^*\mathcal{H} \rightarrow 0$$
is exact, and thus the short exact sequence of chain complexes 
$$ 0 \rightarrow   \widetilde{C}^\bullet({\U, \mathcal{F}}) \rightarrow \widetilde{C}^\bullet({\U, \mathcal{G}}) \rightarrow \widetilde{C}^\bullet({\U, \mathcal{H}}) \rightarrow 0$$
is exact by Corollary \ref{wniosek z surjekcja na sekcjach}. Therefore the proposition follows directly.

\end{proof}
\end{prop}

Hence and by Corollaries \ref{cor - pd F =0 pdG <2} and \ref{wniosek z iniektywnoscia pullbacku dla snopow}, we immediately obtain

\begin{cor}\label{cor- long exact sequence of cohomology}

Let $\U$ be a finite Zariski open covering of $X$ and
$$0 \rightarrow \F \rightarrow \mathcal{G} \rightarrow \mathcal{H} \rightarrow 0$$
be a short exact sequence of quasi-coherent sheaves on $X$. Assume that one of the following conditions holds
\begin{itemize}
\item[i)]
$\pd \F_x=0$ and $\pd \G_x \leq 1$ for all $x \in X$.
\item[ii)]
$\hdim \F =0$ and $\hdim \G \leq 1$.
\end{itemize}

Then there is an induced long exact sequence of blown-up \v{C}ech cohomology with respect to $\U$

$$ \dots \rightarrow \widetilde{H}^p(\U, \F) \rightarrow \widetilde{H}^p(\U, \G) \rightarrow \widetilde{H}^p(\U, \HH) \rightarrow \widetilde{H}^{p+1}(\U, \F) \rightarrow $$
$$ \widetilde{H}^{p+1}(\U, \G) \rightarrow \widetilde{H}^{p+1}(\U, \HH) \rightarrow \widetilde{H}^{p+2}(\U, \F) \rightarrow \dots $$
\hfill$\square$
\end{cor}

\vspace{1ex}

We now prove a version of Cartan's Theorem B for blown-up \v{C}ech cohomologies for some quasi-coherent subsheaves of $\OO_X^{\oplus I}$.

\begin{prop}\label{Cartan B dla pokrycia}
Let $\F$ be a quasi-coherent subsheaf of $ \OO_X^{\oplus I}$ with all stalks $\F_x$, $x\in X$ being free over $\OO_{x,X}$. For any finite Zariski open covering of $X$ we have
$$\widetilde{H}^q(\U, \F)=0.$$

\begin{proof}

We have a short exact sequence of quasi-coherent sheaves
$$ 0 \rightarrow \F \rightarrow \OO_X^{\oplus I} \rightarrow \OO_X^{\oplus I} / \F \rightarrow 0.$$
By assumptions and Corollary \ref{cor - pd F =0 pdG <2}, for any multi-blowup $\alpha : X_\alpha \rightarrow X$ we have a short exact sequence

$$0 \rightarrow \alpha^*\F \rightarrow  \OO_{X_\alpha}^{\oplus I} \rightarrow \OO_{X_\alpha}^{\oplus I}/\alpha^*\F=\alpha^*(\OO_X^{\oplus I} / \F) \rightarrow 0,$$
hence $\alpha^*\F \subset \OO_{X_\alpha}^{\oplus I}$.

Any $[f] \in \widetilde{C}^q(\U, \F)$, has a representative  $f \in C^q(\U^\alpha, \alpha^*)$ for some multi-blowup $\alpha : X_\alpha \rightarrow X$. Our objective is to find a multi-blowup $\beta : X_\beta \rightarrow X_\alpha$ and a $(q-1)$-cocycle $k $ such that $dk=\beta^*f$. To simplify the proof we  assume that $f \in C^q(\U, \F)$. Each of the sets $U_i$ is of the form $U_i = X\setminus \{Q_i=0 \}$ for some $Q_i \in \OO_X(X)$. Recall that $U_{i_0 i_1 \dots i_q}=U_{i_1}\cap U_{i_2}\cap \dots \cap U_{i_q}$ and $f=(f_{i_0 i_1 \dots i_q})$ over all $(q+1)$-tuples of indices for which $i_0,i_1,\dots, i_q \in \{1,2, \dots, n \}$.

By the assumptions every  $f_{i_0 i_1 \dots i_q}$ 
can be identified with a finite tuple of non-vanishing regular functions on $U_{i_0 i_1 \dots i_q}$ i.e
$$f_{i_0 i_1 \dots i_q} \in \F(U_{i_0 i_1 \dots i_q}) \subset \OO_X^{\oplus I}(U_{i_0 i_1 \dots i_q})$$
where
 $$U_{i_0 i_1 \dots i_q}=X\setminus\{Q_{i_0 i_1 \dots i_q}=0 \}$$
  and
   $$Q_{i_0 i_1 \dots i_q}=Q_{i_0} Q_{i_1} \dots Q_{i_q}.$$
   
    By Lemma \ref{q-coh ext section to global} there exist a multi-blowup $\alpha : X_\alpha \rightarrow X$, a positive integer $N_1$ and a global section 
 $g_{i_0 i_1 \dots i_q} \in \alpha^*\OO^{\oplus I}_{X_\alpha}(X_\alpha)$ such that for every $(q+1)$-tuple
$$g_{i_0 i_1 \dots i_q}|_{U^\alpha_{i_0 i_1 \dots i_q}}= (Q^{N_1}_{i_0 i_1 \dots i_q})^\alpha \alpha^*f_{i_0 i_1 \dots i_q}.$$
Consider the image of $g_{i_0 i_1 \dots i_q}$ in the quasi-coherent sheaf $\OO^{\oplus I}_{X_\alpha}/\alpha^*\F$. It is a global section which vanishes on $U^\alpha_{i_0 i_1 \dots i_q}$ and thus, by Lemma \ref{q-coh ext zero section}, there exist a multi-blowup $\beta : X_\beta \rightarrow X_\alpha$ and a positive integer $M$ such that 
$$ (Q_{i_0 i_1 \dots i_q}^M)^{\alpha \circ \beta}  \beta^*g_{i_0 i_1 \dots i_q}$$
is a zero section in $(\beta^*(\OO^I_{X_\alpha}/\alpha^*\F))(X_\beta).$
 Setting $N=N_1+M$, we get 
$$h_{i_0 i_1 \dots i_q}:= (Q_{i_0 i_1 \dots i_q}^M)^{ \alpha \circ \beta}  \beta^* g_{i_0 i_1 \dots i_q}=$$
 $$(Q_{i_0 i_1 \dots i_q}^N)^{\alpha \circ \beta}  {(\alpha \circ \beta)}^*f_{i_0 i_1 \dots i_q} \in (\alpha \circ \beta)^*\F(X_\beta).$$

Since these sections are global, one can always increase $N$ and assume that the number $N$ is even. Define the global regular function 
$$R:= \frac{1}{\sum_{i=1}^n(Q_i^N)^{\alpha \circ \beta}}.$$
Now we are able to define $k \in C^{q-1}(\U^{\alpha \circ \beta}, ({\alpha \circ \beta})^*\F)$ by the formula 
$$ k_{i_0 i_1 \dots i_{q-1}} = R \sum_{i=1}^n\frac{h_{ i i_0 i_1 \dots i_{q-1}}}{(Q^N_{i_0 i_1 \dots i_{q-1}})^{\alpha \circ \beta}}\big|_{U^{\alpha \circ \beta}_{ i_0 i_1 \dots i_{q-1}}}.$$
By the very definition of the operator $d$ we have
$$(dk)_{i_0 i_1 \dots i_q} =\sum_{j=0}^q (-1)^j R \sum_{i=1}^n \frac{h_{ i i_0 i_1 \dots \widehat{i_j} \dots i_{q-1}}}{(Q^N_{i_0 i_1 \dots \widehat{i_j} \dots i_{q-1}})^\beta}\big|_{U^\beta_{ i_0 i_1 \dots i_q}} $$
right hand side is a  finite tuple of regular functions on ${U^{\alpha \circ \beta}_{ i_0 i_1 \dots i_q}}$.
To finish the proof we have to show that $dk=({\alpha \circ \beta})^*f$. It is enough to show it on $\bigcap \U^{\alpha \circ \beta}=U$, since we are dealing only with rational functions and $U$ is Zariski open and dense in ${U^{\alpha \circ \beta}_{ i_0 i_1 \dots i_q}}$.
Recall that, since $({\alpha \circ \beta})^*f$ is a cocycle, we have 
$$0 =(d({\alpha \circ \beta})^*f)_{ i i_0 i_1 \dots i_q} =  ({\alpha \circ \beta})^*f_{i_0 i_1 \dots i_q} + \sum_{j=0}^q (-1)^{j+1} ({\alpha \circ \beta})^*f_{i i_0 i_1\dots \widehat{i_j} \dots i_q}. $$

The following equality holds on $U$:
$$ k_{i_0 i_1 \dots i_{q-1}}= R \sum_{i=1}^n (Q_i^N)^{\alpha \circ \beta} ({\alpha \circ \beta})^*f_{ i i_0 i_1 \dots i_{q-1}}, $$
whence

\begin{equation}\label{equation 1} (dk)_{i_0 i_1 \dots i_q} = \sum_{j=0}^q (-1)^j R \sum_{i=1}^n (Q_i^N)^{\alpha \circ \beta}  ({\alpha \circ \beta})^*f_{ i i_0 i_1\dots \widehat{i_j} \dots i_q}
\end{equation}
Combining the equality (\ref{equation 1}) with the fact that $({\alpha \circ \beta})^*f$ is a cocycle, we obtain
$$(dk)_{i_0 i_1 \dots i_q}= R \sum_{i=1}^n (Q_i^N)^\beta \sum_{j=0}^q (-1)^j ({\alpha \circ \beta})^*f_{ i i_0 i_1\dots \widehat{i_j} \dots i_q}=$$ 
 $$= R \sum_{i=1}^n (Q_i^N)^{\alpha \circ \beta} ({\alpha \circ \beta})^*f_{i_0 i_1 \dots i_q}=({\alpha \circ \beta})^*f_{i_0 i_1 \dots i_q};$$
the last equality follows from the definition of the function $R$. This completes the proof.

\end{proof}
\end{prop}

Now we can readily prove the following real algebraic version of Cartan's Theorem B.

\begin{thm}\label{Cartan B dla pokrycia i hdim =1}
Let $\U=\{U_i \}_{i=1}^n$ be a finite Zariski open covering of $X$ and let $\HH$ be a quasi-coherent sheaf which admits a global presentation such that $\pd \HH_x \leq 1$ for all $x \in X$. Then  $\widetilde{H}^q(\U, \HH) =0$ for $q>0$.

\begin{proof}

We have a short exact sequence
$$ 0 \rightarrow \F \rightarrow \OO_X^{\oplus I} \rightarrow \HH \rightarrow 0. $$
Obviously, $\hdim (\OO_{X}^{\oplus I}) =0$ and  $\pd\HH_x \leq 1$ by assumption. Hence $\pd \F_x=0$ for all $x \in X$ by  Proposition \ref{Prop z 4 podpunktami dla modulow}. By Corollary \ref{cor- long exact sequence of cohomology}, the above exact sequence induces the long exact sequence of blown-up \v{C}ech cohomology
$$ \dots \rightarrow \widetilde{H}^q(\U, \F) \rightarrow \widetilde{H}^q(\U, \OO_{X}^{ \oplus I}) \rightarrow \widetilde{H}^q(\U, \HH) \rightarrow \widetilde{H}^{q+1}(\U, \F)$$
$$ \rightarrow \widetilde{H}^{q+1}(\U, \OO_{X}^{ \oplus I}) \rightarrow \widetilde{H}^{q+1}(\U, \HH) \rightarrow \dots $$
By Proposition \ref{Cartan B dla pokrycia} 
$$\widetilde{H}^q(\U, \F) = \widetilde{H}^q(\U, \OO_{X}^{ \oplus I})= \widetilde{H}^{q+1}(\U, \F) = \widetilde{H}^{q+1}(\U, \OO_{X}^{ \oplus I}) =0$$
for any $q>0$. Hence $\widetilde{H}^q(\U, \HH)=0$ as asserted.

\end{proof}
\end{thm}

\begin{cor}\label{cartan b dla hdim=1 coh}
Let $\U=\{U_i \}_{i=1}^n$ be a finite Zariski open covering of $X$ and let $\HH$ be a coherent sheaf of homological dimension $\leq 1$. Then  $\widetilde{H}^q(\U, \HH) =0$ for $q>0$.

\begin{proof}
By Lemma \ref{Lemma 8} there exists a multi-blowup $\alpha : X_\alpha \rightarrow X$ such that $\alpha^*\HH$ admits a global presentation. We can thus apply Theorem \ref{Cartan B dla pokrycia i hdim =1}

\end{proof}

\end{cor}

The family of finite Zariski open coverings of $X$ can be directed by refinement relation $\succeq$. Let $\U=\{U_i\}_{i=1}^n$ be the refinement of $\V=\{V_j\}^m_{j=1}$ i.e. $\U \succeq \V$. Let $\tau : \{1,2, \dots, m \} \rightarrow \{1,2,\dots,n \}$ be the refinement map. Then the  homorphism
$$ \widetilde{C}^\bullet (\V, \F) \rightarrow \widetilde{C}^\bullet (\U, \F)$$
does depend on the choice of $\tau$.
However, by \cite[I, \S3, Prop.\ 3]{Serre} the induced map on the cohomology
$$\widetilde{H}^q(\V,\F) \rightarrow \widetilde{H}^q(\U,\F)$$
does not depend on the choice of $\tau$. Hence we may define

\emph{The $q$-th blown-up \v{C}ech cohomology group  $\widetilde{H}^q(X, \F)$ of $\F$} is the direct limit
$$\widetilde{H}^q(X, \F)=\lim_{ \underset{\U}{\rightarrow}} \widetilde{H}^q(\U,\F).$$

Summing up, Theorem \ref{Cartan B dla pokrycia i hdim =1}, Proposition \ref{Cartan B dla pokrycia} and Corollary \ref{cartan b dla hdim=1 coh} yield immediately the following general version of Cartan's Theorem B.

\begin{thm}
Let $\F$ be a sheaf of $\OO_X$-modules and let $\U$ be a finite Zariski open covering of $X$. Assume that one of the following conditions hold
\begin{itemize}
\item[a)]
$\F$ is a quasi-coherent subsheaf of $ \OO_X^{\oplus I}$ such that $\pd \F_x =0$ for all $x\in X$.

\item[b)]
$\F$ is a quasi-coherent sheaf of global presentation such that $\pd \F_x \leq 1$ for all $x \in X$.

\item[c)]
$\F$ is a coherent sheaf and $\hdim \F \leq 1$.
\end{itemize}
Then $\widetilde{H}^q(\U, \F)=0$ and, a fortiori, $\widetilde{H}^q(X, \F)=0$ for $q \geq 1$.

\hfill$\square$
\end{thm}

\section{A real algebraic version of the first Cousin problem}

In this section we deal with the first Cousin problem. The classical complex version of the first Cousin problem is treated e.g.\ in \cite{grauert}. Before discussing details, we give an outline of the problem.
 Let $\U = \{U_i \}_{i=1}^n$ be a finite Zariski open covering of $X$. Assume that for each $i$ we have a rational function $f_i$ on $U_i$ such that $f_i - f_j $ is regular on $U_i \cap U_j$ for each two distinct indices $i,j=1,2, \dots n$. Then we call $\{(U_i,f_i) \}_{i=1}^n$ data of the first Cousin problem or an additive Cousin distribution on $X$.

Let $U$ be a Zariski open subset of $X$. We say that two rational functions $f,g$ on $X$ have the same \emph{principal part} on $U$ if $f-g \in \OO_X(U)$.

The first Cousin problem consists in characterizing those data
\newline
 $ \{ (U_i,f_i)\}_{i=1}^n$ which have the principal parts of a rational function $f$ on $X$, i.e.\ those for which $f-f_i$ are regular on $U_i$, $i=1,2,\dots, n$. We then say that the data 
$\{(U_i,f_i) \}_{i=1}^n$ is \emph{solvable}. If every first Cousin data on $X$ is solvable, we say that the first Cousin problem is \emph{universally solvable on X} (see e.g.\ \cite{grauert}, Introduction, \S2).

We are going to describe the above problem in terms of sheaves. Let $\KK_X$ be the constant sheaf of rational functions on $X$. Consider the short exact sequence   
 $$ 0 \rightarrow \OO_X \rightarrow \KK_X \overset{\varphi}{\rightarrow} \HH_X :=\KK_X/\OO_X \rightarrow 0 $$
of quasi-coherent sheaves on $X$. The data $\{(U_i,f_i) \}_{i=1}^n$ of the first Cousin problem can be related to a unique global section $s \in \HH_X(X)$; every such section is called a principal part distribution on $X$. Then we also say that 
$\{(U_i,f_i) \}_{i=1}^n$ is an $s$-representing distribution. In particular, for every rational function $f \in \KK_X(X)$ we have its principal part distribution $\varphi(f)$ on $X$. 

Let $\{(U_i,f_i) \}_{i=1}^n$ be an $s$-representing distribution. Then any rational function $f \in \KK_X(X)$ satisfying $\varphi(f)=s$ is one such that $f-f_i \in \OO_X(U_i)$ for each $i=1,2,\dots,n$.

Take the short exact sequence of quasi-coherent sheaves
$$ 0 \rightarrow \OO_X \rightarrow \KK_X \rightarrow \HH_X \rightarrow 0. $$
Then we have merely the following exact sequence of chain complexes
$$ 0 \rightarrow C^\bullet (\U, \OO_X) \rightarrow C^\bullet(\U, \KK_X) \rightarrow C^\bullet(\U, \HH_X).$$
Under the circumstances, we have the following exact sequence for the classical \v{C}ech cohomologies (cf.\ \cite[Chap I, \S 3, section 24]{Serre})
$$  0 \rightarrow H^0(X,\OO_X) \rightarrow H^0(X,\KK_X) \rightarrow H^0(X,\HH_X) \overset{\zeta}{\rightarrow} $$
$$ H^1(X,\OO_X) \rightarrow H^1(X,\KK_X) \rightarrow H^1(X,\HH_X). $$ 
The above sequence of chain complexes induces a long exact sequence for \v{C}ech cohomologies whenever the topology of $X$ is paracompact (cf.\ \cite[Chap I, \S 3, section 25]{Serre}). This can be applied in the classical case of complex analytic geometry (see e.g.\ \cite{grauert}).

Clearly, every $s$-representing distribution $\{(U_i,f_i) \}_{i=1}^n$ determines a $1$-cocycle $(g_{ij})$, $g_{ij}:=f_i-f_j$, which induces a cohomology class $\zeta(s)$ in $H^1(\U, \OO_X)$. The solvability of a given Cousin data can be rephrased in terms of vanishing $\zeta(s)$ in $H^1(\U, \OO_X)$.
\begin{lemma}
An $s$-representing distribution $\{(U_i,f_i) \}_{i=1}^n$ is solvable iff $\zeta(s)=0 \in H^1(\U, \OO_X)$.

\begin{proof}
The necessary condition.
Let $f$ be a rational function such that $f-f_i \in \OO_X(U_i)$. Then $(g_{ij})$ is a coboundary of a $0$-cocycle $h$ given by $h_i=f_i-f$, hence $\zeta(s)=0$.

Conversely, if $\zeta(s)$=0 then $g=dh$ for some $0$-cocycle $(h)_i$. We then have $(h_i-h_j)|_{U_i \cap U_j}=(f_i-f_j)|_{U_i \cap U_j}$, and obviously $(f_i-h_i)|_{U_i\cap U_j}=(f_j-h_j)|_{U_i\cap U_j}$ for any pair of two indices $i,j$. It follows that the system $(f_i-h_i)$ for $i=1,2,\dots,n$ can be glued to a global rational function $f$ such that $f-f_i \in \OO_X(U_i)$, as asserted.

\end{proof}

\end{lemma}

 No natural map from $H^q(X,\OO_X)$ to $H^q(\U, \OO_X)$ exist for arbitrary covering $\U$. Therefore it is necessary to refine a given covering to solve the first Cousin problem in complex analytic geometry (cf.\ \cite{grauert}, Chapter V). In view of Theorem 5.2., it turns out to be superfluous in real algebraic geometry after blowing up. This is stated in the following main

\begin{thm}
Let $\{(U_i,f_i) \}_{i=1}^n$ be an $s$-representing distribution. Then there exists a multi-blowup $
\alpha : X_\alpha \rightarrow X$ such that pull-back $\{(U^\alpha_i,f^\alpha_i) \}_{i=1}^n$ is solvable i.e.\ there exists a rational function $f$ on $X_\alpha$ such that 
\newline 
${f-f_i^\alpha \in \OO_{X_\alpha}(U_i^\alpha)}$.

\end{thm}

We need an elementary lemma

\begin{lemma}\label{inj in frac field}
For any multi-blowup $\alpha:X_\alpha \rightarrow X$, we have a short exact sequence 
$$0 \rightarrow \OO_{X_\alpha} \rightarrow \KK_{X_\alpha} \rightarrow \HH_{X_\alpha} \rightarrow 0$$
of quasi-coherent sheaves on $X_\alpha$.

\begin{proof}
For any $y \in X_\alpha, x=\alpha(y)\in X$ the inclusion $\OO_{x,X} \hookrightarrow \KK_x=K$ induces the inclusion
$$  \OO_{x,X} \otimes_{\OO_{x,X}} \OO_{y,X_\alpha}=\OO_{y,X_\alpha} \hookrightarrow \KK_x\otimes_{\OO_{x,X}} \OO_{y,X_\alpha}=K=\KK_y. $$
Hence the conclusion follows.

\end{proof}
\end{lemma}

\noindent
{\em Proof of Theorem 6.1} 
 The above lemma yields the following short exact sequence of chain complexes
 $$ 0 \rightarrow \widetilde{C}^\bullet(\U, \OO_X) \rightarrow \widetilde{C}^\bullet(\U, \KK_X) \rightarrow  \widetilde{C}^\bullet(\U, \HH_X) \rightarrow 0$$
 which induces a long exact sequence in blown-up \v{C}ech cohomology
 $$ 0 \rightarrow \widetilde{H}^0(\U, \OO_X) \rightarrow \widetilde{H}^0(\U, \KK_X) \rightarrow \widetilde{H}^0(\U, \HH_X) \rightarrow \widetilde{H}^1(\U, \OO_X)$$
  $$ \rightarrow \widetilde{H}^1(\U, \KK_X) \rightarrow \widetilde{H}^1 (\U, \HH_X) \rightarrow \dots $$
 
 By Proposition \ref{Cartan B dla pokrycia},  $\widetilde{H}^1(\U, \OO_X)=0$. Since 
 $$ \widetilde{H}^1(\U, \OO_X) = \LIM H^1(\U^\alpha, \OO_{X_\alpha}),$$ 
 for any class $\omega \in H^1(\U, \OO_X)$ there exists a multi-blowup $\alpha : X_\alpha \rightarrow X$ such that $\alpha^*\omega=0 $ in $ H^1(\U^\alpha\, \OO_{X_\alpha})$. To finish the proof it is enough to take $\omega = \zeta(s)$.

\hspace*{\fill} $\square$
\newline
 
\begin{remark}
It follows from the above proof that there is an isomorphism of $\OO_X(X)$-modules
$$\widetilde{H}^q(\U, \KK_X) \cong \widetilde{H}^q(\U, \HH_X)$$
for $q \geq 1$.

\end{remark}

 \section{An Example}
 In this section we provide an example (based on \cite[Example 12.1.5]{Bochnak}) of a quasi-coherent locally free sheaf $\F$ on $\mathbb{R}^2$ of infinite rank such that for any multi-blowup $\alpha : X_\alpha \rightarrow \mathbb{R}^2$ the pull-back sheaf $\alpha^*\F$ is not generated by global sections. This shows that the version of Cartan's Theorem A from our paper \cite{TK} cannot be generalized to quasi-coherent sheaves. 
 
For $k,l \in \mathbb{N}$, the irreducible polynomial 
$$P_{k,l}(x,y)=x^{2k}(x-1)^{2l} +y^2$$
has at most two zeros, $c_1=(0,0)$ and $c_2=(1,0)$.  Put ${U_i := \mathbb{R}^2\setminus \{ c_i\}},\newline { i=1,2}$. Then the transition function
$$g_{2,1}: U_1\cap U_2 \rightarrow \mathrm{GL}(1,\mathbb{R})=\mathbb{R}^*$$
$$(x,y) \mapsto P_{k,l}(x,y)$$
determines a vector bundle $\xi_{k,l}$ of rank $1$ on $\mathbb{R}^2$. The sheaf $\F_{k,l}$ of its sections is locally free of rank 1. A global section $s$ of $\xi_{k,l}$ can be described as a pair of regular functions
$$ s_i : U_i \rightarrow \mathbb{R}, \,\,\, i=1,2$$
such that $g_{2,1}s_1=s_2$. It is clear that the bundles $\xi_{0,l}$ and $\xi_{k,0}$ are trivial line bundles.

Suppose now that $k,l >0$ and set $s_i=\frac{f_i}{h_i}$ where $f_i$ and $h_i$ are relatively prime polynomials. Then $P_{k,l}f_1h_2=f_2h_1$. Since $P_{k,l}$ cannot divide either $h_1$ or $h_2$, we get $f_2=\lambda P_{k,l} f_1$ and $h_1=\lambda^{-1}h_2$ with some $\lambda \in \mathbb{R}, \lambda \neq 0$. Therefore any global section $s$ has to vanish at $c_1=(0,0)$, and thus $\F_{k,l}$ cannot be generated by global sections.

Similarly, if we put $g_{1,2}=P_{k,l}(x,y)$ in the above construction, then we obtain a locally free sheaf of rank one such that every global section vanishes at $c_2$.

The case $(k,l)=(1,1)$ is treated in \cite{TK}. It is showed there that if $\sigma : X_\sigma \rightarrow \mathbb{R}^2$ is a blowup of the plane at the origin, then the global section of $\sigma^*\F_{1,1}$ given by the pair $\widetilde{s}=( \widetilde{s_1},\widetilde{s_2})$ with 
$$\widetilde{s_1}=\frac{1}{(x^2+y^2)^\sigma} \ \ \text{and} \ \  \widetilde{s_2}=\frac{P_{1,1}^\sigma(x,y)}{(x^2+y^2)^\sigma}$$
does not vanish on $X_\sigma$.

 Let $\tau : X_\tau \rightarrow \mathbb{R}^2$ be a blowup of the plane at the point $(1,0)$. In a similar way, the global section $\widetilde{t}=(\widetilde{t_1}, \widetilde{t_2})$ given by 
 $$\widetilde{t_1}=\frac{((x-1)^2 +y^2)^\tau}{P_{1,1}^\tau(x,y)} \  \text{and} \ \ \widetilde{t_2}=((x-1)^2 +y^2)^\tau$$
  is a nowhere vanishing global section of $\tau^*\F_{1,1}$.

\begin{lemma}\label{lemma 7.1}
Let $k,l \in \mathbb{N}, k,l>0$, $d_1,d_2,\dots,d_r$ be a finite number of points of $\mathbb{R}^2$ distinct from $c_1$ and $c_2$, $\sigma : X_\sigma \rightarrow \mathbb{R}^2$ be a multi-blowup at $d_1,d_2,\dots, d_r$ and $D=\{d_1,d_2, \dots, d_r \}$. Then 
$$\sigma^*\xi_{k,l}|_{X_\sigma \setminus \sigma^{-1}(D)} \cong \xi_{k,l}|_{\mathbb{R}^2\setminus D}$$
and the pull-back $\sigma^*\xi_{k,l}$ is not generated by global sections.

\begin{proof}
Since a blowup is biregular off the exceptional divisor, it is enough to show that $\xi_{k,l}|_{\mathbb{R}^2\setminus D}$ is not generated by global sections. Let $V_i:=U_i\setminus D$, $i=1,2$. Then the reasoning from the beginning of this section can be repeated verbatim with the $U_i$ replaced by the $V_i$, $i=1,2$.
\end{proof}
\end{lemma}

Before stating the next lemma, we recall the construction of the blowup at the origin $\sigma : X_\sigma \rightarrow \mathbb{R}^2$.
We have 
$$X_\sigma = \{(x,y, u:v) \in \mathbb{R}^2 \times \mathbb{P}^1(\mathbb{R}) : xv=uy \}.$$
We can cover $X_\sigma$ with two Zariski open subsets
 $$\Omega_1= \{(x,y,u:v)\in X_\sigma :u \neq 0 \}\,\,\, \text{and} \,\,\, \Omega_2 = \{(x,y,u:v) \in X_\sigma : v \neq0\}$$
 with local coordinates $$(x,\frac{v}{u}) \ \ \text{and} \ \ (\frac{u}{v},y)$$ respectively. In these local coordinates $\sigma$ is expressed by the formulas
$$\sigma(r,s)=(r,rs) \ \ \text{and} \ \ \sigma(r,s)=(rs,s),$$
respectively.

\begin{lemma}\label{lemma 7.2}
The pull-back $\sigma^*\xi_{k,l}$ is a trivial line bundle on $\Omega_2$ and $\sigma^*\xi_{k,l}|_{\Omega_1}$ is isomorphic to $\xi_{k-1,l}$ on $\Omega_1 \cong \mathbb{R}^2_{r,s}$.

\begin{proof}
The first assertion is obvious since $c_2 \notin \sigma(\Omega_2)$. To verify the second one, compute the transition function of $\sigma^*\xi_{k,l}|_{\Omega_1}$:
$$\widetilde{g_{21}}(r,s)=P_{k,l}^\sigma(r,s) = r^2(r^{2k-2}(r-1)^{2l} +s^2)=r^2P_{k-1,l}(r,s)$$
on $\Omega_1$. But the line bundle $\eta$ on $\Omega_1$ with the transition function $r^2$ is trivial on $\Omega_1$. Therefore
$$\sigma^*\xi_{k,l}|_{\Omega_1} \cong \xi_{k-1,l} \otimes \eta \cong \xi_{k-1,l}$$
as asserted.

\end{proof}
\end{lemma}

Similarly, let $\tau : X_\tau \rightarrow \mathbb{R}^2$ be a blowup at the point $c_2$.
We have 
$$X_\tau = \{(x,y, (u-1):v) \in \mathbb{R}^2 \times \mathbb{P}^1(\mathbb{R}) : (x-1)v=(u-1)y \}.$$
We can cover $X_\tau$ by two Zariski open sets 
$$\Omega_1= \{(x,y,(u-1):v)\in X_\tau : u\neq 1 \}\,\, \text{and} \,\,\Omega_2 = \{ (x,y,(u-1):v)\in X_\tau : v\neq 0 \}$$
 with local coordinates
$$(x,\frac{v}{u-1}) \,\,\, \text{and} \,\,\, (\frac{u-1}{v},y)$$
respectively. In these local coordinates $\tau$ is expressed by the formulas 
$$\tau(r,s)=(r,(r-1)s)\,\,  \text{and} \,\, \tau(r,s)=(rs+1,s)$$ respectively. As before, we obtain

\begin{lemma}\label{lemma 7.3}

The pull-back $\tau^*\xi_{k,l}$ is a trivial line bundle on $\Omega_2$ and $\tau^*\xi_{k,l}|_{\Omega_1}$ is isomorphic to $\xi_{k,l-1}$ on $\Omega_1 \cong \mathbb{R}^2_{r,s}$.

\end{lemma}
\hfill$\square$

In view of Lemmas \ref{lemma 7.1}, \ref{lemma 7.2}, \ref{lemma 7.3} it is clear that blowing up at a point $d\neq c_1,c_2$ is immaterial and what improves the bundle $\xi_{k,l}$ is only successive blowing up at the points $(0,0)$ or $(1,0)$ on the chart $\Omega_1$. Each such blowup transforms the initial line bundle (isomorphic to $\xi_{p,q}$) to the line bundle $\xi_{p-1,q}$ or $\xi_{p,q-1}$, respectively, on the chart $\Omega_1$. We must continue until we attain the line bundle $\xi_{0,q}$ or $\xi_{p,0}$, with some $0 \leq p \leq k, 0 \leq q \leq l$, which is generated by global sections. In this manner, we have proven the following

\begin{prop}

Let $\sigma : X_\sigma \rightarrow \mathbb{R}^2$ be a composition of $r$ blowups and $k,l \in \mathbb{N}$. If the pull-back $\sigma^*\xi_{k,l}$ or, equivalently, $\sigma^*\F_{k,l}$ is generated by global sections, then $r \geq \min (k,l)$.

\end{prop}
\hfill$\square$

\begin{remark}
The composition $\alpha$ or $\beta$ of $k$ or $l$ blowups $\sigma$ or $\tau$ at the points $(0,0)$ or $(1,0)$, respectively, on the successive charts $\Omega_1$ transforms $\xi_{k,l}$ to a trivial line bundle. Indeed, it is not difficult to check that the nowhere vanishing section on $X_\alpha$ is given by 
$\widetilde{s}=(\widetilde{s_1},\widetilde{s_2})$

$$\widetilde{s_1}=\frac{1}{(x^{2k}+y^2)^\alpha} \, \, \text{and} \,\, \widetilde{s_2}=\frac{P^\alpha_{k,l}}{(x^{2k}+y^2)^\alpha}.$$
And that the global nowhere vanishing section on $X_\beta$ is given by $\widetilde{t}=(\widetilde{t_1},\widetilde{t_2})$ 
$$\widetilde{t_1}=\frac{((x-1)^{2l} +y^2)^\beta}{P_{k,l}^\beta(x,y)} \  \text{and} \ \ \widetilde{t_2}=((x-1)^{2l} +y^2)^\beta.$$

\end{remark}

We immediately obtain

\begin{cor}
Let 
$$\F := \bigoplus_{k=1}^\infty \F_{k,k}$$
be a quasi-coherent locally free sheaf of infinite rank on $\mathbb{R}^2$. Then for any multi-blowup $\alpha : X_\alpha \rightarrow \mathbb{R}^2$ the pull-back $\alpha^*\F$ is not generated by global sections.

\end{cor}

 \section*{Acknowledgement}
 The author would like to kindly thank Professor Krzysztof Nowak for his useful comments and discussions.

\vspace{1ex}

\begin{small}
\noindent
Tomasz Kowalczyk

\noindent
Institute of Mathematics

\noindent
Faculty of Mathematics and Computer Science

\noindent
Jagiellonian University

\noindent
ul. Profesora Łojasiewicza 6

\noindent
 30-348 Kraków, Poland

\noindent
e-mail: tomek.kowalczyk@student.uj.edu.pl

\end{small}

\end{document}